\documentclass[twoside]{article}

\usepackage{amsmath,amssymb,amsfonts}
\usepackage{amsthm}
\usepackage{xspace,enumerate}
\usepackage{burdges}
\usepackage{calc}

\newcommand\Uir{U_{0,r}}

\begin{document}

\title{The Bender method in groups\\ of finite Morley rank}
\author{Jeffrey Burdges
\thanks{Supported by an NSF Graduate Research Fellowship,
 NSF grant DMS-0100794, and
 Deutsche Forschungsgemeinschaft grant Te 242/3-1.}\\
Fakult\"at f\"ur Mathematik \\ Universit\"at Bielefeld \\ D-33501 Bielefeld, Germany \\
{\tt burdges@math.rutgers.edu}
}

\markboth%
{\hfill J.~Burdges \hfill}%
{\hfill The Bender method \hfill}
\pagestyle{myheadings}


\maketitle

The algebraicity conjecture for simple groups of finite Morley rank,
also known as the Cherlin-Zilber conjecture, states that simple groups
of finite Morley rank are simple algebraic groups over algebraically
closed fields.  In the last 15 years, the main line of attack on this problem
has been Borovik's program of transferring methods from finite group
theory.  Borovik's program has led to considerable progress; however,
the conjecture itself remains decidedly open.
In Borovik's program, groups of finite Morley rank are divided into
four types, odd, even, mixed, and degenerate, according to the structure
of their Sylow 2-subgroup.  For {\em even} and {\em mixed type}
the algebraicity conjecture has been proven.
The present paper is part of the program to bound the {\em \Prufer rank}
of minimal simple groups of finite Morley rank and {\em odd type}.

In \cite{CJ01}, Cherlin and Jaligot achieved a bound of \Prufer rank two
for {\em tame} minimal simple groups.  Here a group of finite Morley rank
is said to be tame if it does not involve a field of finite Morley rank with
a proper infinite definable subgroup of it's multiplicative group.
Cherlin, Jaligot, and the present author will bound the \Prufer rank
at two in \cite{BCJ}.

Tameness is used in two important ways in \cite{CJ01}.  The final
number theoretic contradiction of \cite{CJ01} uses tameness in an
essential way, and \cite{BCJ} will completely replace this argument.
However, the very first use of tameness in \cite{CJ01} produces the
following fact, which shows that intersections of Borel subgroups are abelian.

\begin{namedtheorem}{Jaligot's Lemma}[{\cite[\qLemma 3.11]{CJ01}}]%
\label{named:Ojaligot}
Let $G$ be a {\em tame} minimal connected simple group of finite Morley rank.
Let $B_1$ and $B_2$ be two distinct Borel subgroups of $G$ with
$O(B_1) \neq 1$ and $O(B_2) \neq 1$.  Then $F(B_1) \cap F(B_2) = 1$.
\end{namedtheorem}

The present paper examines the worst violations of Jaligot's Lemma
in the {\em nontame} context, i.e.\ those involving nonabelian
intersections of Borel subgroups.  We fail to exclude such nonabelian
intersections outright, but we gain much information about the
specific local configuration responsible for nonabelian intersections.

In the context of minimal simple groups, the present paper provides a
analog of Bender's Uniqueness Theorem \cite[\qTheorem 28.2]{GLS1}
(see also \cite[Ch.\ 5]{Gagen} and \cite[\S9]{BenderGlauberman}),
a result underlying the Bender method \cite[\S28]{GLS1} of analyzing
the maximal subgroups containing the centralizer of an involution.  
Both the Bender Uniqueness Theorem and the present paper provide
information about the normalizers of various subgroups of the intersection
of two distinct maximal subgroups.
However, our situation will be simplified by two facts:
torsion behaves extremely well (see \S\ref{sec:Upjaligot}), and
our ``torsion-free primes'', so-called reduced ranks, are naturally
ordered by their degree of unipotence (see Fact \ref{nilpotencepre2}).

In \cite{BCJ}, much of the information about this nonabelian configuration,
plus analysis of the relevant abelian intersections, is used to prove the following.

\begin{namedtheorem}{Theorem}
Let $G$ be a minimal connected simple group of finite Morley rank
 and odd type with a strongly embedded subgroup.
Then $G$ has \Prufer rank one.
\end{namedtheorem}

\noindent  One proves the bound on \Prufer rank by showing that
simple groups of finite Morley rank and \Prufer rank at least three
have strongly embedded subgroups \cite{BuPhd}.

The bulk of this paper consists of the analysis of nonabelian
maximal intersections of Borel subgroups in a minimal connected
simple group of finite Morley rank (see \S\ref{sec:Aj}).
A priori, the analysis of these maximal intersections yields only
a limited description of nonmaximal intersections.

\begin{namedtheorem}{Proposition \ref{AjaligotC}}
Let $G$ be a minimal connected simple group of finite Morley rank.
Let $B_1,B_2$ be two distinct Borel subgroups of $G$, and
 $H$ a connected subgroup of the intersection $B_1 \cap B_2$.
Then the following hold.
\begin{conclusions}
\item $H'$ is rank homogeneous for $r' := \rr_0(H')$.
\item Every connected nilpotent subgroup of $H$ is abelian.
\item $F_{r'}(H) = U_{0,r'}(H)$ is a Sylow $U_{0,r'}$-subgroup of $H$.
\item[(etc.)]
\end{conclusions}
\end{namedtheorem}

In high \Prufer rank, the experience of \cite{BCJ} suggests that
Proposition \ref{AjaligotC} itself is insufficient, but that the results of
\S\ref{sec:Aj} which describe the configuration arising from maximal
nonabelian intersections {\em are} sufficient.
In practice, these results can be used because of the following equivalence
between different characterizations of nonabelian maximal intersections.

\begin{namedtheorem}{Theorem \ref{maximalities}}
Let $G$ be a minimal connected simple group of finite Morley rank,
 and let $B_1,B_2$ be two distinct Borel subgroups of $G$.
Suppose that $H := (B_1 \cap B_2)^\o$ is nonabelian.
Then the following are equivalent.
\begin{conclusions}
\item $B_1$ and $B_2$ are the only Borel subgroups of $G$ containing $H$.
\item If $B_3$ and $B_4$ are distinct Borel subgroups of $G$ containing $H$,
  then\\ $(B_3 \cap B_4)^\o = H$.
\item If $B_3 \neq B_1$ is a Borel subgroup containing $H$,
  then $(B_3 \cap B_4)^\o = H$.
\item $C^\o_G(H')$ is contained in $B_1$ or $B_2$.
\item $B_1$ and $B_2$ are not conjugate under $C^\o_G(H')$.
\item $\rr_0(B_1) \neq \rr_0(B_2)$.
\end{conclusions}
\end{namedtheorem}

The results presented here do not assume the presence of 2-torsion.  As such,
we expect these results to play a significant role in the study of both odd and
degenerate type groups.

The paper begins by recalling the necessary background in
\S\ref{sec:Background}, including the definition of 0-unipotence.
Section \ref{sec:Upjaligot} proves Jaligot's Lemma for Borel
subgroups with $p$-unipotent radicals, and thus eliminates most
concerns with connected torsion.  Section \ref{sec:Aj} carries out the core
of our analysis of a maximal nonabelian intersection of Borel subgroups.
Section \ref{sec:Conclusions} proves the equivalence of the various
notions of nonabelian maximal intersection (see Theorem \ref{maximalities}),
and summarizes the results of Section \ref{sec:Aj} in that context
(see Theorem \ref{summary}).
Section \ref{sec:Future} discusses possible some future directions
related to Carter subgroups.

\section{Background}\label{sec:Background}

\subsection{Unipotent groups}

While there is no intrinsic definition of unipotence in
 a group of finite Morley rank,
there are various analogs of the ``unipotent radical'':
the Fitting subgroup, the $p$-unipotent operators $U_p$, for $p$ prime,
 and their ``charateristic zero'' analogs $\Uir$ from \cite{Bu03,BuPhd}.
We recall their definitions.

\begin{definition}
The {\em Fitting subgroup} $F(G)$ of a group $G$
 of finite Morley rank is the subgroup generated
 by all its nilpotent normal subgroups.
\end{definition}

The Fitting subgroup is itself nilpotent and definable
 \cite[Theorem 7.3]{Bel87,Ne91,BN},
and serves as a notion of unipotence in some contexts.
However, since the Fitting subgroup of a solvable group may not be
contained in the Fitting subgroups of a solvable group containing it,
it is not a robust notion.

\begin{fact}[{\cite[\qCorollary 9.9, \qTheorem 9.21]{BN}}]\label{fittingquotient}
Let $G$ be a connected solvable group of finite Morley rank.
Then $G' \leq F^\o(G)$ and $G/F^\o(G)$ is a divisible abelian group.
\end{fact}

\begin{fact}[{\cite[\qLemma 4.20]{Fre00b}; compare \cite[Ex. 3 p.\ 148]{BN}}]%
\label{cent_divtorsion}\label{cent_fittingtorsion}
Let $G$ be a connected solvable group of finite Morley rank,
and let $T$ be a divisible torsion subgroup of $G$.
Then $T \cap F(G)$ is central in $G$ 
\end{fact}

\begin{definition}
A subgroup of a connected solvable group $H$
 of finite Morley rank is said to be {\it $p$-unipotent}
if it is a definable connected $p$-group of bounded exponent.
\end{definition}

\begin{fact}
[{\cite[\qCorollary 2.16]{CJ01}; \cite[\qFact 2.36]{ABC97}}]
\label{Upnilpotence} 
Let $H$ be a connected solvable group of finite Morley rank.
Then there is a unique maximal $p$-unipotent subgroup of $H$,
denoted $U_p(H)$, and $U_p(H) \leq F^\o(H)$.
\end{fact}

The $p$-unipotent radical $U_p$ will automatically behave well under
intersections with other solvable groups.

\begin{fact}[{\cite[\qTheorem 9.29 and \S6.4]{BN}}]\label{Sylow_conplus}
Let $G$ be a connected solvable group of finite Morley rank. 
Then a Sylow $p$-subgroup $P$ of $G$ is connected, and
 $P = U_p(G) * T$ for a divisible abelian $p$-group $T$.
\end{fact}


%

The present paper relies on the theory of ``characteristic zero''
unipotence introduced in \cite{Bu03}.
We now turn our attention to this (long) definition,
 as well as some facts from \cite{Bu03,Bu05a,BuPhd}.

We say that a connected abelian group of finite Morley rank is
{\it indecomposable} if it has a unique maximal proper
definable connected subgroup, denoted $J(A)$.

\begin{fact}[{\cite[\qLemma 2.4]{Bu03}}]\label{decomposition}
Every connected abelian group of finite Morley rank can be written
as a finite sum of definable indecomposable abelian subgroups.
\end{fact}

\begin{definition}
We define the {\it reduced rank} $\rr(A)$ of a definable abelian
group $A$ to be the Morley rank of the quotient $A/J(A)$,
i.e.\ $\rr(A) = \rk(A/J(A))$.
For a group $G$ of finite Morley rank, and any integer $r$, we define
$$ \Uir(G) = \Genst{A \leq G}{%
\parbox{\widthof{$A$ is a definable indecomposable group,}}%
{$A$ is a definable indecomposable group, \\
\hspace*{10pt} $\rr(A) = r$, and $A/J(A)$ is torsion-free}}\mathperiod $$
We say $G$ is a {\it $\Uir$-group} (alternatively {\it $(0,r)$-unipotent})
if $U_{0,r}(G)=G$.
The {\em 0-unipotent radical} $U_0(G)$ is the nontrivial $U_{0,r}(G)$ with
$r$ maximal.  We also set $\rr_0(G) = \max \{r \mid \Uir(G) \neq 1 \}$.
So $U_0(G) = U_{0,\rr_0(G)}(G)$.

As a notational convention, we define $F_r(G) = U_{0,r}(F(G))$, and use
$\Ftor(G)$ to denote the definable closure of the torsion subgroup of $F(G)$.
\end{definition}

We view the reduced rank parameter $r$ as a {\em scale of unipotence},
 with larger values being more unipotent.
By the following fact, the ``most unipotent'' groups, in this scale, are nilpotent.

\begin{fact}%
[{\cite[\qTheorem 2.21]{BuPhd}; \cite[\qTheorem 2.16]{Bu03}}]
\label{nilpotence}
Let $H$ be a connected solvable group of finite Morley rank.
Then $U_0(H) \leq F(H)$.
\end{fact}

\begin{fact}[{\cite[\qLemma 2.5]{Bu05a}}]\label{nilpotencecor}
Let $H$ be a connected solvable group of finite Morley rank.
Let $r$ be the maximal reduced rank such that $\Uir(H) \not\leq Z_n(H)$
for all $n$.  Then $\Uir(H) \leq F(H)$.
\end{fact} 

The two preceding facts prove nilpotence of $\Uir(H)$ for the largest
values of $r$.  However, $H$ may intersect a solvable group $K$ with
$\rr_0(K) > \rr_0(H)$.  Many other facts about $\Uir$-groups hold
for arbitrary values of $r$

\begin{fact}%
[{\cite[\qLemma 2.12]{BuPhd}; \cite[\qLemma 2.11]{Bu03}}]
\label{Uhom} 
Let $f : G \to H$ be a definable homomorphism between two groups of finite
Morley rank.  Then the following hold.
\begin{conclusions}
\item \textit{(Push-forward)}
  $f(\Uir(G)) \leq \Uir(H)$ is a $\Uir$-group.
\item \textit{(Pull-back)}
  If $\Uir(H) \leq f(G)$ then $f(\Uir(G)) = \Uir(H)$.
\end{conclusions}
In particular, an extension of a $\Uir$-group by a $\Uir$-group is a
$\Uir$-group.
\end{fact}

\begin{fact}%
[{\cite[\qLemma 2.26]{BuPhd}; compare \cite[\qLemma 2.3]{Bu05a}}]
\label{Ucenter}
Let $G$ be a nilpotent group of finite Morley rank satisfying
$\Uir(G) \neq 1$ or $U_p(G) \neq 1$.  Then
$\Uir(Z(G)) \neq 1$ or $U_p(Z(G)) \neq 1$, respectively.
\end{fact}

We have a 0-unipotent analog of the connected normalizer condition
of \cite[Lemma 6.3]{BN}.

\begin{fact}%
[{\cite[\qLemma 2.28]{BuPhd}; \cite[\qLemma 2.4]{Bu05a}}]
\label{Unormalizer} 
Let $G$ be a nilpotent $\Uir$-group of finite Morley rank.
If $H < G$ is a definable subgroup then $\Uir(N_G(H)/H) > 1$.
\end{fact}

Our next result generalizes the fact that a finite nilpotent group is the product
of its Sylow $p$-subgroups.

\begin{fact}%
[{\cite[\qCorollary 3.6]{Bu05a}; \cite[\qTheorem 2.31]{BuPhd}}]
\label{nildecomp}
Let $G$ be a connected nilpotent group of finite Morley rank.
Then $G = D * B$ is a central product of definable characteristic
 subgroups $D,B \leq G$ where $D$ is divisible and
 $B$ is connected of bounded exponent.
Let $T$ be the torsion part of $D$.
Then we have decompositions of $D$ and $B$ as follows.
\begin{eqnarray*}
D &=& d(T) * U_{0,1}(G) * U_{0,2}(G) * \cdots \\
B &=& U_2(G) \times U_3(G) \times U_5(G) \times \cdots
\end{eqnarray*}
\end{fact}

Here $d(T)$ denotes the definable closure of $T$, which is defined
 to be the intersection of all definable subgroups containing $T$.

\begin{fact}%
[{\cite[\qLemma 2.32]{BuPhd}; \cite[\qCorollary 3.7]{Bu05a}}]
\label{Unilcommutator}
Let $G$ be a solvable group of finite Morley rank, let $S \subseteq G$ be
any subset, and let $H$ be a nilpotent $\Uir$-group which is normal in $G$.
Then $[H,S] \leq H$ is a $\Uir$-group.
\end{fact}

Olivier \Frecon has improved this result in \cite{Fre06}.

Our next fact says that ``more unipotent'' groups do not act
 on ``less unipotent'' groups.

\begin{fact}%
[{\cite[\qLemma 4.4]{Bu05a}; see also \cite[\qCorollary 3.8]{FrJa04}}]
\label{nilpotencepre2}
Let $G = H T$ be a group of finite Morley rank with
 $H \normal G$ a nilpotent $U_{0,r}$-group and
  $T$ a nilpotent $U_{0,s}$-group for some $s \geq r$.
Then $G$ is nilpotent.
\end{fact}

In \cite{Wag01}, Wagner showed that fields of finite Morley rank and
characteristic $p \neq 0$ have no torsion free sections of their multiplicative
groups \cite{Wag01}.  The Zilber Field Theorem \cite[Theorem 9.1]{BN}
allows us to rephrase Wagner's result as follows.

\begin{fact}[{\cite[\qLemma 4.3]{Bu05a}}]\label{Uz_on_Up}
Let $G$ be a connected solvable group of finite Morley rank.
Suppose that $S$ is a nilpotent $\Uir$-subgroup of $G$, and
 that $G = U_p(G) S$ for some $p$ prime.
Then $G$ is nilpotent, and $[U_p(G),S]=1$.
\end{fact}


\subsection{Toral groups}

\begin{definition}
A definable subgroup $C$ of a group $G$ of finite Morley rank
 which is nilpotent and almost self-normalizing in $G$ is called
 a {\em Carter subgroup} of $G$.
\end{definition}

The following result is a summary, in order, of \cite[Proposition 3.2]{Fre00a},
\cite[Corollary 4.8]{Fre00a}, \cite[Theorem 5.5.12]{Wag,Wag94},
and \cite[Corollary 7.15]{Fre00a}.

\begin{fact}
\label{Carter_cleanup}\label{Carter_con}\label{Carter_exists}
\label{Carter_Sylow}\label{Carter_conj}\label{Carter_cover}
Let $H$ be a connected solvable group of finite Morley rank. Then the following hold.
\begin{conclusions}
\item $H$ has a Carter subgroup.
\item The Carter subgroups of $H$ are the definable nilpotent
subgroup of $H$ with $N_H^\o(Q) = Q$.
In particular, Carter subgroups of $H$ are connected.
\item The Carter subgroups of $H$ are conjugate in $H$.
\item The Carter subgroups of $H$ cover $H/H'$.
\item Let $R$ be a Sylow $p$-subgroup of $H$.
Then $N_H(R)$ contains a Carter subgroup of $H$.
\end{conclusions}
\end{fact}

\begin{definition}
A {\em Sylow $\Uir$-subgroup} of a group $G$ of finite Morley rank
is a maximal definable nilpotent $\Uir$-subgroup of $G$.
\end{definition}

Sylow $\Uir$-subgroups are an analog of Carter subgroups in the following sense.

\begin{lemma}[{\cite[\qLemma 4.18]{BuPhd}; \cite[\qLemma 5.2]{Bu05a}}]
\label{USylow}
Let $H$ be a group of finite Morley rank.
Then the Sylow $\Uir$-subgroups of $H$ are exactly those
 nilpotent $\Uir$-subgroups $S$ such that $\Uir(N_G(S)) = S$.
\end{lemma}

\begin{fact}[{\cite[\qLemma 4.19]{BuPhd}; \cite[\qTheorem 5.7]{Bu05a}}]\label{USylow_decomp}
Let $H$ be a connected solvable group of finite Morley rank and
let $Q$ be a Carter subgroup of $H$.
Then $\Uir(H') \Uir(Q)$ is a Sylow $\Uir$-subgroup of $H$, and
every Sylow $\Uir$-subgroup has this form for some Carter subgroup $Q$.
\end{fact}

\begin{fact}[{\cite[\qTheorem 4.16]{BuPhd}; \cite[\qTheorem 5.5]{Bu05a}}]\label{UirCarter_conj}
Let $H$ be a connected solvable group of finite Morley rank.
Then the Sylow $\Uir$-subgroups of $H$ are $H$-conjugate.
\end{fact}

\section{$p$-Unipotence}\label{sec:Upjaligot}

In this section, we show that intersections of the Fitting subgroups
of distinct Borel subgroups are torsion free, and thus eliminate
many concerns about torsion from the main analysis to follow.
The arguments of this section are based directly on the original
proof of Jaligot's lemma \cite[\qLemma 3.11]{CJ01} (see introduction).

\begin{lemma}\label{Upjaligot}
Let $G$ be a minimal connected simple group of finite Morley rank.
Let $B_1,B_2$ be two distinct Borel subgroups satisfying $U_{p_i}(B_i) \neq 1$
for some prime $p_i$ ($i=1,2$).  Then $F(B_1) \cap F(B_2) = 1$.
\end{lemma}

\begin{proof}
We first show that $U_p(B_1 \cap B_2) = 1$ for all $p$, prime.
Suppose toward a contradiction that $X := U_p(B_1 \cap B_2) \neq 1$.
We may assume that $\rk(X)$ is maximal among all choices of $B_1$ and
$B_2$.  Let $B$ be a Borel subgroup of $G$ containing $N^\o_G(X)$.

We now show that $B_i = B$ for $i=1,2$.
If $X = U_p(B_i)$, then $B \geq N^\o_G(U_p(B_i)) = B_i$, and $B = B_i$.
So we consider the case where $X < U_p(B_i)$. 
By Fact \ref{Upnilpotence}, the group $U_p(B_1)$ is nilpotent.
By the normalizer condition \cite[Lemma 6.3]{BN},
$$ X < U_p(N_{U_p(B_i)}^\o(X)) \leq U_p(B)\mathperiod $$
Since $B$ is a Borel subgroup containing $N_{U_p(B_i)}^\o(X)$,
 $B_i = B$ by the maximality of $\rk(X)$.
Thus $B_1 = B = B_2$, a contradiction.
So $U_p(B_1 \cap B_2) = 1$.

We now prove the lemma.  Suppose toward a contradiction that there
is an $f \in F(B_1) \cap F(B_2)$ with $f \neq 1$.  By Fact \ref{Ucenter},
$Z_i := U_{p_i}(Z(F(B_i)))$ is a nontrivial $U_{p_i}$-subgroup
of $Z(F(B_i))$ for $i=1,2$.  Let $B$ be a Borel subgroup containing $C^\o_G(f)$.
As $Z_i \leq Z(F(B_i)) \leq C^\o_G(f)$ for $i=1,2$, we find
$Z_i \leq U_{p_i}(B_i \cap B)$ for $i=1,2$.
Thus $B_1 = B = B_2$ by the first part.
\end{proof}

The above conclusion holds, with a similar proof, if we replace $U_p$ by
$U_{0,\rr_0(G)}$; however, one does not know that all Borel subgroups satisfy
 $\rr_0(B) = \rr_0(G)$, just as one does not know that $U_p(B) \neq 1$.
These techniques can be extended to eliminate all
 torsion from the intersection of the Fitting subgroups.

\begin{corollary}\label{Xtf}
Let $G$ be a minimal connected simple group of finite Morley rank.
Let $B_1,B_2$ be two distinct Borel subgroups of $G$.
Then $F(B_1) \cap F(B_2)$ is torsion free.
\end{corollary}

\begin{proof}
Suppose towards a contradiction that, for some prime $p$,
 there is a nontrivial Sylow $p$-subgroup $P$ of $X := F(B_1) \cap F(B_2)$.
We may assume $U_p(B_2) = 1$ by Lemma \ref{Upjaligot}.
So $P$ is central in $B_2$ by Fact \ref{cent_fittingtorsion},
and hence $C^\o_G(P) = B_2$.
If $U_p(B_1) = 1$ too, then $P$ is central in $B_1$
 by Fact \ref{cent_fittingtorsion}.  Thus $U_p(B_1) \neq 1$.
Let $R$ be a Sylow $p$-subgroup of $B_1$ containing $P$.
By Fact \ref{Sylow_conplus}, $R = U_p(B_1) * T$
 for some divisible abelian $p$-group $T$.
By Fact \ref{Ucenter}, $Y := U_p(C_{U_p(B_1)}(P)) \neq 1$.
But $Y \leq B_2$, a contradiction.
\end{proof}

\section{Maximal intersections}\label{sec:Aj}

In this section, we analyze intersections of Borel subgroups
which are maximal in the sense of the following definition.

\begin{definition}
Let $G$ be a minimal connected simple group of finite Morley rank,
 and let $B_1,B_2$ be two distinct Borel subgroups of $G$.
We say that $H := (B_1 \cap B_2)^\o$ is a {\em maximal
intersection} if $H$ is maximal among all choices of distinct
Borel subgroups $B_1$ and $B_2$.  In this situation, we refer
to $B_1,B_2$ as a {\em maximal pair} of Borel subgroups of $G$.
\end{definition}

The analysis of maximal intersections will not always directly produce
information about general intersections of Borel subgroups.
One hopes to produce results which translate down, in some form, to
nonmaximal intersections.  For example, our lives would be simple if
the intersection $H$ turned out to be abelian.
However, such a simple analog of Jaligot's Lemma eludes us.
Instead, we discover a plausible nonabelian maximal configuration,
which this section explores in detail.

In the specific case of \cite{BCJ}, the configuration below survives
until the end of the analysis, and eventually dies for the same reasons as
the case of an abelian intersection does.  Since roughly half of the facts
below are important in \cite{BCJ}, we view them collectively as a
machine for handling nonabelian intersections.

Throughout this section, we consider a minimal connected simple
group $G$ of finite Morley rank, and a maximal pair $B_1,B_2$
of Borel subgroups of $G$ which violate Jaligot's Lemma.

\begin{hypothesis}\label{Aj_hypothesis}
We assume the following.
\begin{hypotheses}
\item  $G$ is a minimal connected simple group of finite Morley rank.
\item  $B_1$ and $B_2$ are a maximal pair of Borel subgroups of $G$.
\item  The intersection $F(B_1) \cap F(B_2)$ of
   their Fitting subgroups is nontrivial.
\end{hypotheses}
\end{hypothesis}

\begin{notation}\label{Aj_notation}
We also adopt the following notation.
\begin{hypotheses}
\item $H := (B_1 \cap B_2)^\o$ denotes the maximal intersection.
\item $X := F(B_1) \cap F(B_2)$ denotes the intersection of the Fitting subgroups.
\item $r' := \rr_0(X)$ denotes the reduced rank of $X$.
\end{hypotheses}
\end{notation}

In addition, we frequently discuss the reduced rank $r'$ piece $Y := F_{r'}(H)$
of the Fitting subgroup of $H$.
A Carter subgroup $Q$ of $H$ will also play a central role.

We observe that $\rr_0(H') = r'$ if $H' \neq 1$ by Theorem \ref{Aj_UrX} below,
 which explains our choice of notation.

\subsection{Examples}

Before beginning our analysis in \S\ref{subsec:Rank_homogeneity},
we describe a few ``nearly algebraic'' configurations which survive.
This material will not be used below, but may shed some light on
our goals.

We consider two Borel subgroups $B_1$ and $B_2$ such that
 $H := (B_1 \cap B_2)^\o$ is a maximal intersection.
We suppose that $\rr_0(B_1) \geq \rr_0(B_2)$, and that
$H$ is nonabelian.  The ``light'' Borel subgroup $B_2$ could be the
group of upper triangular $3 \times 3$ matrices over an
algebraically closed field $k$ of characteristic zero, and say of determinant one. 
The important requirement here is that $F(B_2)$ is nonabelian.

It will be shown that in fact $\rr_0(B_1) > \rr_0(B_2)$.
After breaking the symmetry in this way, it turns out that there are
striking differences between the ``heavy'' $B_1$ and the ``light'' $B_2$.
In particular, the Borel subgroup $B_2$ can be algebraic over an
algebraically closed field $k$ of characteristic zero, while $B_1$
necessarily interprets a bad field.

One way to proceed is as follows.  Let $B_2$ be the subgroup of
upper triangular $3 \times 3$ matrices with $a_{11} = a_{33} = 1$.
$$ B_2 = \begin{pmatrix} 1 & k_+ & k_+ \\ & T & k_+ \\ & & 1 \\ \end{pmatrix} $$
So $B_2 = U \rtimes T$ where $U$ is the (unipotent) group of
strictly upper triangular matrices, and $T$ is a one dimensional torus.

Before specifying $B_1$, we need to choose $H = (B_1 \cap B_2)^\o$
so that $H = N_G(H')$ and $F(H)$ is abelian.  A suitable choice of
$H$ is the subgroup of $B_2$ given by $a_{23} = 0$.
Then $H = F(H) \rtimes T$ with $F(H) = H' \times Z(U)$.
$$ H = \begin{pmatrix} 1 & H' & Z(U) \\ & T & 0 \\ & & 1 \\ \end{pmatrix} $$

One possibility for $B_1$ would be the direct product
$$ B_1 = [H' \rtimes T] \times [U_0(B_1) \rtimes Z(U)] $$
where $H' \rtimes T$ behaves as it does in $B_2$, while $Z(U)$
has become part of a bad field, and $U_0(B_1)$ is an additive group
of larger reduced rank.

In the above situation, a Carter subgroup $Q = T \times Z(U)$
of $H$ is a Carter subgroup of both $B_1$ and $B_2$.
We will see below that $Q$ will always be a Carter subgroup of $B_1$.
However, a Carter subgroup of $H$ need not be a Carter subgroup
of $B_2$, in general.

For example, we may take $B_2$ to be the Borel subgroup of $\SL_3(k)$.
Here we must take $H'$ to be a one dimensional unipotent subgroup of $B_2$
which is normalized by some one dimensional torus $T$, but not by the full torus.
As above, $H = N_{B_2}(H')$ and $T \cong (H' \rtimes T) \times Z(U)$.
Indeed, $Q = T \times Z(U)$ is a Carter subgroup of $H$, but $Z(U)$ is
no longer part of a Carter subgroup of $B_2$, and $T$ is a proper subgroup
of a Carter subgroup of $B_2$.  As a consequence $T$ may not centralize
$U_0(B_1)$.  Indeed, the Carter subgroup $Q = T \times Z(U)$ of $H$
(and $B_1$) may be the full multiplicative group of a bad field, while only
$T$ acts on $H'$!


\subsection{Rank homogeneity of $X$}\label{subsec:Rank_homogeneity}

We have two goals in the first stage of the analysis.
First, we will show that $\rr_0(B_1) \neq \rr_0(B_2)$,
a key fact in the remainder of the analysis.
Second, we will show that the subgroup $X := F(B_1) \cap F(B_2)$
is rank homogeneous in the sense of the following definition.

\begin{definition}
A group $K$ of finite Morley rank is said to be {\em rank homogeneous}
if $K$ is torsion free and $U_{0,r}(K) = 1$ for $r \neq \rr_0(K)$.
\end{definition}

We observe that $X \normal H$, and that $H' \leq X$.
So rank homogeneity of $X$ will imply $r' = \rr_0(H')$ too, if $H' \neq 1$.
Clearly $r' > 0$ by Corollary \ref{Xtf}.  In particular, $\rr_0(H) > 0$

To begin our analysis, we may assume that
$$ \rr_0(B_1) \geq \rr_0(B_2)\mathperiod \eqno{\star} $$

As a first step, the normalizer condition shows that the reduced rank must
grow on one side of our intersection. 

\begin{lemma}\label{Aj_B1}
$\rr_0(H) < \rr_0(B_1)$
\end{lemma}

\begin{proof}
Suppose toward a contradiction that $\rr_0(H) = \rr_0(B_1)$.
Since $\rr_0(H) \leq \rr_0(B_2) \leq \rr_0(B_1)$ by our assumption,
all are equal and $U_0(H) \leq U_0(B_1) \cap U_0(B_2)$.
Let $B_3$ be a Borel subgroup of $G$ containing $N^\o_G(U_0(H))$.

We now show that $B_i = B_3$ for $i=1,2$.
If $U_0(H) = U_0(B_i)$, then $B_3 \geq N^\o_G(U_0(H)) = B_i$,
 and $B_3 = B_i$.
So we suppose that $U_0(H) < U_0(B_i)$.
By Fact \ref{nilpotence}, $U_0(B_i)$ is nilpotent.
By Fact \ref{Unormalizer},
 $$ U_0(H) < U_{0,\rr_0(H)}(N_{U_0(B_i)}(U_0(H)) \leq B_3\mathperiod $$
Since $U_0(H) \normal H$, $(B_i \cap B_3)^\o > H$.
By maximality of $H$, $B_i=B_3$ here too.
Thus $B_1=B_3=B_2$, a contradiction.
\end{proof}

Starting with our next lemma, we use the decomposition of nilpotent groups,
given in Fact \ref{nildecomp}, to ``blow up'' the centralizers of various
subgroups of $H$.  This is a variation on the normalizer condition based
argument used above.

Some of our next lemmas use the decomposition of nilpotent groups
(Fact \ref{nildecomp}) instead of the connected normalizer condition
of \cite[Lemma 6.3]{BN}.
Such arguments resemble the use of Fact \ref{cent_divtorsion} in Corollary \ref{Xtf}.

\begin{lemma}\label{Aj_B2}
$\rr_0(H) = \rr_0(B_2)$
\end{lemma}

\begin{proof}
Suppose toward a contradiction that $\rr_0(H) < \rr_0(B_2)$.
Since $U_0(X) \neq 1$ by Corollary \ref{Xtf},
 there is a Borel subgroup $B_3$ of $G$ containing $N^\o_G(U_0(X))$.
Since $U_0(X) \normaleq H$, we have $H \leq B_3$.
Since $\rr_0(H) < \rr_0(B_i)$ for $i=1,2$ by Lemma \ref{Aj_B1}
 and our hypothesis, we have
 $U_0(B_i) \leq C^\o_G(U_0(X)) \leq B_3$ by Fact \ref{nildecomp}.
So $(B_i \cap B_3)^\o > H$ for $i=1,2$,
 and hence $B_1=B_3=B_2$, a contradiction.
\end{proof}

So $\rr_0(B_1) > \rr_0(B_2)$ by Lemmas \ref{Aj_B1} and \ref{Aj_B2}.
In the event that $H$ is nonabelian, this rank inequality prevents other
 Borel subgroups from containing $H$.

\begin{proposition}\label{Aj_B1B2unique}
If $H$ is nonabelian, then
$B_1$ and $B_2$ are the only Borel subgroups containing $H$.
\end{proposition}

\begin{proof}
Suppose toward a contradiction that there is a Borel subgroup $B_3$,
distinct from $B_1$ and $B_2$, which contains $H$.
By the maximality of $H$,
 $(B_1 \cap B_3)^\o = H$ and $(B_2 \cap B_3)^\o = H$.
Since $H' \leq F(B_3)$,  the maximal pairs
 $B_1,B_3$ and $B_3,B_2$ satisfy Hypothesis \ref{Aj_hypothesis}.
Since $\rr_0(B_1) > \rr_0(H)$ by Lemma \ref{Aj_B1},
 $\rr_0(B_3) = \rr_0(H)$ by Lemma \ref{Aj_B2}.
But, since $\rr_0(B_2) = \rr_0(H)$ by Lemma \ref{Aj_B2},
 $\rr_0(B_3) > \rr_0(H)$ by Lemma \ref{Aj_B1},
a contradiction.
\end{proof}

The groups $B_1$ and $B_2$ are not conjugate, since
 $\rr_0(B_1) > \rr_0(B_2)$, a point exploited by the following lemma.

\begin{lemma}\label{Aj_FBi}
$F^\o(B_i) \not\leq H$ for $i=1,2$.
\end{lemma}

\begin{proof}
Since $\rr_0(H) < \rr_0(B_1)$ by Lemma \ref{Aj_B1},
 we have $F^\o(B_1) \not\leq H$.
Suppose toward a contradiction that $F^\o(B_2) \leq H$.
Then $H \normal B_2$ by Fact \ref{fittingquotient}.
So $H \leq B_1 \cap B_1^g$ for $g\in B_2 \setminus B_1$,
contradicting Lemmas \ref{Aj_B1} and \ref{Aj_B2}.
\end{proof}
%

Our analysis hinges upon understanding the behavior of the normalizers
of various $U_{0,r'}$-subgroups of $H$.

\begin{lemma}\label{Aj_NX}
For any nontrivial definable $X_1 \leq X$ with $X_1 \normal H$,
 we have $N^\o_G(X_1) \leq B_1$.
\end{lemma}

\begin{proof}
Let $B_3$ be a Borel subgroup containing $N^\o_G(X_1)$.
Then $H \leq B_3$.
Since $\rr_0(H) < \rr_0(B_1)$ by Lemma \ref{Aj_B1},
Fact \ref{nildecomp} yields
 $U_0(B_1) \leq C^\o_G(X_1) \leq B_3$;
 and thus $(B_1 \cap B_3)^\o > H$.
By the maximality of $H$,  $B_1=B_3$.
So $N^\o(X_1) \leq B_1$.
\end{proof}

In particular, the previous two lemmas show that $X \cap Z(F(B_2)) = 1$.

Our results now diverge from the conclusions of
 Bender's Uniqueness Theorem \cite[\qTheorem 28.2]{GLS1}
in that only one ``prime'' may occur in $X$.

\begin{theorem}\label{Aj_UrX}
$X$ is rank-homogeneous.
In particular, $X = U_0(X)$.
\end{theorem}

\begin{proof}
By Corollary \ref{Xtf}, $X$ is torsion free.
Suppose toward a contradiction that $U_{0,r}(X) \neq 1$
 for some $r < r'$.
By Fact \ref{nildecomp} and Lemma \ref{Aj_NX},
$$ F(B_2) \leq C^\o(U_{0,r'}(X)) C^\o(U_{0,r}(X))
   \leq B_1 $$
and $F(B_2) \leq H$, contradicting Lemma \ref{Aj_FBi}.
Hence $X = U_0(X)$.
\end{proof}

\subsection{Fitting subgroup of $B_2$}

Our next goal is to understand the Fitting subgroup of  $B_2$.
In particular, we will determine which factors of $F(B_2)$ are
contained in $H$.

\begin{lemma}\label{Aj_B2div}
$F^\o(B_2)$ is divisible, and $\Ftor^\o(B_2) \leq Z(H)$.
\end{lemma}

\begin{proof}
Of course, $U_p(B_2) \leq \Ftor^\o(B_2)$ by Fact \ref{Upnilpotence}.
By Theorem \ref{Aj_UrX} and Fact \ref{nildecomp},
 $\Ftor^\o(B_2) \leq N^\o(X)$.
Since $N^\o(X) \leq B_1$ by Lemma \ref{Aj_NX},
 $\Ftor^\o(B_2) \leq H$.
So $U_p(B_2) \leq H$ and $\Ftor^\o(B_2) \leq F^\o(H)$.
But $U_p(H) = 1$ for all primes $p$, by Lemma \ref{Upjaligot}.
By Fact \ref{Sylow_conplus}, $F^\o(B_2)$ is divisible.
So $\Ftor^\o(B_2)$ is central in $H$ by Fact \ref{cent_divtorsion}.
\end{proof}

\begin{lemma}\label{Aj_UrFB2}
$F_r(B_2) \leq Z(H)$ for $r\neq r'$.
\end{lemma}

\begin{proof}
By Fact \ref{nildecomp} and Theorem \ref{Aj_UrX},
 $F_r(B_2) \leq C^\o_G(X)$.
Since $C^\o_G(X) \leq B_1$ by Lemma \ref{Aj_NX},
 $F_r(B_2) \leq H$.
Since $F_r(B_2)$ is normalized by $H \leq B_2$,
we have $F_r(B_2) \leq F(H)$.
By Fact \ref{Unilcommutator} and Theorem \ref{Aj_UrX},
\[ [H,F_r(B_2)] \leq U_{0,r}(H') \leq U_{0,r}(X) = 1 \qedhere \]
\end{proof}

We know that one part of the Fitting subgroup of $B_2$
 is not contained in $H$.

\begin{lemma}\label{Aj_UFB2}
$F_{r'}(B_2) \not\leq H$ is not abelian.
\end{lemma}

\begin{proof}
By Lemma \ref{Aj_FBi}, $F^\o(B_2) \not\leq H$.
Since $\Ftor^\o(B_2) \leq H$ by Lemma \ref{Aj_B2div},
 $F_s(B_2) \not\leq H$ for some $s$ by Fact \ref{nildecomp}.
So $F_{r'}(B_2) \not\leq H$ by Lemma \ref{Aj_UrFB2}.
Since $N^\o_G(X) \leq B_1$ by Lemma \ref{Aj_NX},
 $F_{r'}(B_2)$ is not abelian.
\end{proof}

Lemmas \ref{Aj_UrFB2} and \ref{Aj_UFB2}
tell us that $r'$ is uniquely determined by $B_2$.

\begin{corollary}\label{Aj_B2rprime}
$F_r(B_2)$ is nonabelian iff $r = r'$.
\end{corollary}

\begin{proof}
Immediate from Lemmas \ref{Aj_UrFB2} and \ref{Aj_UFB2}.
\end{proof}

The first two lemmas of this section give us a measure of control over
the large reduced ranks in $B_2$.

\begin{lemma}\label{Aj_UrB2cor}
$U_{0,r}(B_2) \leq F(B_2)$ for $r > r'$.
\end{lemma}

\begin{proof}
Let $A$ be an indecomposable $U_{0,r}$-subgroup of $B_2$.
Since $r > r'$, the group
 $A \cdot F_{r'}(B_2)$ is nilpotent by Fact \ref{nilpotencepre2}.
By Fact \ref{nildecomp}, $A$ centralizes $F_{r'}(B_2)$, including $X$.
So $A \leq H$ by Lemma \ref{Aj_NX}.
By Lemma \ref{Aj_B2div}, $A$ centralizes $\Ftor^\o(B_2)$.
By Lemma \ref{Aj_UrFB2}, $A$ centralizes $F_s(B_2)$ for $s \neq r'$.
So $A$ centralizes $F^\o(B_2)$ by Fact \ref{nildecomp},
 and hence $A \cdot F^\o(B_2)$ is nilpotent.
By Fact \ref{fittingquotient}, $A \cdot F^\o(B_2) \normal B_2$,
 and hence $A \leq F(B_2)$.
\end{proof}

In particular, $F_r(B_2) = U_{0,r}(B_2)$
 is the unique Sylow $U_{0,r}$-subgroup of $H$.
We also observe that, if $H$ contains a Carter subgroup of $B_2$, then
 $r = r'$ is the maximal reduced rank such that $\Uir(H) \not\leq Z(H)$,
and hence $U_{0,r'}(B_2) \leq F(B_2)$ too, by Fact \ref{nilpotencecor}.
One may hope to show that $U_{0,r'}(B_2) \leq F(B_2)$ without assuming
that $H$ contains a Carter subgroup of $B_2$.

\subsection{Structure of $H$}

We now turn our attention towards various subgroups of $H$.
First, if we restrict ourselves to $U_{0,r'}(H)$, the argument of
Lemma \ref{Aj_UrB2cor} applies to the reduced rank $r'$ itself.

\begin{lemma}\label{Aj_Yident}\label{Aj_Yident2}
$U_{0,r'}(H) \leq F(B_2)$.  In particular, $F_{r'}(H) = U_{0,r'}(H)$
 is the unique Sylow $U_{0,r'}$-subgroup of $H$,
\end{lemma}

\begin{proof}
Let $A$ be an indecomposable $U_{0,r'}$-subgroup of $H$.
By Fact \ref{nilpotencepre2}, $A \cdot F_{r'}(B_2)$ is nilpotent.
For $r\neq r'$, $A$ centralizes $F_r(B_2)$,
 by Lemma \ref{Aj_UrFB2}.
By Lemma \ref{Aj_B2div},
 $A$ centralizes $\Ftor^\o(B_2)$ too.
So $A \cdot F^\o(B_2)$ is nilpotent by Fact \ref{nildecomp}.
By Fact \ref{fittingquotient}, $A \cdot F(B_2) \normal B_2$,
 and hence $A \leq F(B_2) \cap H \leq F(H)$.
So $U_{0,r'}(H)\leq F_{r'}(H)$,
 and these two subgroups are equal.
\end{proof}

We adopt the notation $Y := U_{0,r'}(H)$ ($= F_{r'}(H)$).
We find that $Y$ opposes the pull of $X$.

\begin{lemma}\label{Aj_Y}
$N^\o_G(Y) \leq B_2$ and $Y > X$.
In addition, $U_{0,r'}(N_{F(B_2)}(Y)) \not\leq H$.
\end{lemma}

\begin{proof}
Let $P := F_{r'}(B_2)$.
Then $Y \leq P$ by Lemma \ref{Aj_Yident}.
By Lemma \ref{Aj_UFB2}, $F_{r'}(B_2) \not\leq H$,
 so $Y < P$.
By Fact \ref{Unormalizer}, $Y < U_{0,r'}(N^\o_P(Y))$.
Now $Y > X$ by Lemma \ref{Aj_NX},
 and $N^\o_G(Y) \leq B_2$ by maximality of $H$.
\end{proof}  

We can now prove one of our main results.

\begin{theorem}\label{Aj_FH}
Every connected definable nilpotent subgroup of $H$ is abelian.
\end{theorem}

\begin{proof}
Since $U_p(H) = 1$ for any prime $p$ by Lemma \ref{Upjaligot},
 the Sylow $p$-subgroups of $H$ are abelian, by Fact \ref{Sylow_conplus}.
For $r\neq r'$, $U_{0,r}(H') = 1$ by Theorem \ref{Aj_UrX},
 so a Sylow $U_{0,r}$-subgroup of $H$ is abelian by Fact \ref{Unilcommutator}
 (or via Fact \ref{Uhom}).
If $Y' \neq 1$, then
 $N^\o_{B_2}(Y') \leq B_1$ by Lemma \ref{Aj_NX},
 contradicting Lemma \ref{Aj_Y}.
So $Y$ is abelian.
By Lemma \ref{Aj_Yident},
$Y$ is the unique Sylow $U_{0,r'}$-subgroup of $H$.
Thus all Sylow subgroups of $H$ are abelian.
For any connected definable nilpotent subgroup $K$ of $H$,
 $K$ is the central product of its (generalized) Sylow subgroups,
 by Fact \ref{nildecomp}, so $K$ is abelian.
\end{proof}

We observe that $H$ tends to be almost self-normalizing.

\begin{lemma}\label{Aj_NH}
If $H$ is nonabelian, then $N^\o_G(H) = H$.
\end{lemma}

\begin{proof}
By Lemma \ref{Aj_Y}, $N^\o_G(H) \leq N^\o_G(Y) \leq B_2$.
Since $H$ is nonabelian,
 $N^\o_G(H) \leq N^\o_G(H') \leq B_1$ too,
 by Lemma \ref{Aj_NX}.
\end{proof}

\subsection{Structure of $B_1$}

We now turn our attention toward the Fitting subgroup of $B_1$,
 and a Carter subgroup $Q$ of $H$.

\begin{lemma}\label{Aj_B1div}
$F^\o(B_1)$ is divisible, and $\Ftor^\o(B_1) \leq Z(H)$.
\end{lemma}

\begin{proof}
By Fact \ref{Uz_on_Up}, $U_p(B_1)$ is centralized by $Y$.
By Lemma \ref{Aj_Y}, $U_p(B_1) \leq H$.
But $U_p(H) = 1$ for any prime $p$, by Lemma \ref{Upjaligot}. 
By Fact \ref{Sylow_conplus},
 $\Ftor^\o(B_1)$, and $F^\o(B_1)$, are divisible.
So $\Ftor^\o(B_1)$ is central in $B_1$ by Fact \ref{cent_divtorsion}.
By Lemma \ref{Aj_Y}, $\Ftor^\o(B_1) \leq N^\o_G(Y) \leq B_2$,
 and $\Ftor^\o(B_1) \leq Z(H)$.
\end{proof}

We observe that $U_p(B_i) = 1$ for $i=1,2$ and $p$ prime,
 by Lemmas \ref{Aj_B2div} and \ref{Aj_B1div}.

\begin{lemma}\label{Aj_BNX}
$F_{r'}(B_1) = X$.  So $B_1 = N^\o_G(X)$.
\end{lemma}

\begin{proof}
By Fact \ref{nilpotencepre2}, 
 $Y \cdot F_{r'}(B_1)$ is nilpotent.
By Lemma \ref{Aj_Y}, $N^\o_{Y \cdot F_{r'}(B_1)}(Y) \leq H$.
So we obtain $F_{r'}(B_1) \leq Y$ by Fact \ref{Unormalizer}.
By Lemma \ref{Aj_Yident2}, $Y \leq F(B_2)$,
 so $F_{r'}(B_1) \leq X$.
\end{proof}

\begin{corollary}\label{Aj_FB1CX}
$F_{r'}(B_1)$ is abelian, and $F^\o(B_1) \leq C^\o_G(X)$.
\end{corollary}

\begin{proof}
By Theorem \ref{Aj_FH} and Lemma \ref{Aj_BNX},
 $F_{r'}(B_1)$ is abelian, and $F_{r'}(B_1) \leq C^\o_G(X)$.
For $r \neq r'$,
 $F_r(B_1) \leq C^\o_G(X)$ by Fact \ref{nildecomp}.
Also $\Ftor^\o(B_1) \leq C^\o_G(X)$ by Fact \ref{nildecomp}.
So $F(B_1) \leq C^\o_G(X)$ by Fact \ref{nildecomp}.
\end{proof}

We now examine a Carter subgroup $Q$ of $H$.

\begin{lemma}\label{Aj_UQ}
$U_{0,r'}(Q) = U_{0,r'}(Z(H))$, and this group is nontrivial.
\end{lemma}

\begin{proof}
By Lemma \ref{Aj_Y}, $U_{0,r'}(H/H') \neq 1$.
So $U_{0,r'}(Q) \neq 1$ by Facts \ref{Carter_cover}-4 and \ref{Uhom}.
By Theorem \ref{Aj_FH}, $Q$ and $Y$ are abelian.
By Lemma \ref{Aj_Yident}, $U_{0,r'}(Q) \leq Y$.
So $U_{0,r'}(Q)$ centralizes both $Q$ and $H' \leq Y$.
Thus $U_{0,r'}(Q) \leq Z(H)$ by Fact \ref{Carter_cover}-4.
Conversely, $Z^\o(H) \leq Q$.
\end{proof}

\begin{theorem}\label{Aj_Carter}
$N^\o_G(U_{0,r'}(Q)) \leq B_2$.  So $N^\o_G(Q) \leq B_2$,
 and $Q$ is a Carter subgroup of $B_1$.
\end{theorem}

\begin{proof}
We first show that $N^\o_G(U_{0,r'}(Q)) \leq B_2$.
By Lemma \ref{Aj_Y}, $N^\o_G(Y) \leq B_2$.
So we may assume that $U_{0,r'}(Q) < Y$,
 and hence $H$ is nonabelian by Lemma \ref{Aj_UQ}.
So $B_1$ and $B_2$ are the only Borel subgroups of $G$ which contain $H$,
 by Proposition \ref{Aj_B1B2unique}.
By Lemma \ref{Aj_UQ}, $H \leq N^\o_G(U_{0,r'}(Q)) < G$,
 and hence $N^\o_G(U_{0,r'}(Q)) \leq B_1 \cup B_2$.
By Lemma \ref{Aj_BNX}, $N^\o_G(X) = B_1$.
Since $Y = X U_{0,r'}(Q)$ by Fact \ref{USylow_decomp},
 $N^\o_{B_1}(U_{0,r'}(Q)) \leq N^\o_G(Y) \leq B_2$
 (by Lemma \ref{Aj_Y}).
So $N^\o_G(U_{0,r'}(Q)) \leq B_2$.

Thus $N^\o_G(Q) \leq N^\o_G(U_{0,r'}(Q)) \leq B_2$.
By Fact \ref{Carter_cleanup}-2, $Q$ is a Carter subgroup of $B_1$.
\end{proof}

We can now show that $r'$ is the only reduced rank appearing
 in both $F(B_1)$ and $F(B_2)$.

\begin{lemma}\label{Aj_FrB1}
$F_r(B_1) = 1$ for $r \leq \rr_0(B_2)$ with $r \neq r'$.
\end{lemma}

\begin{proof}
Let $T := F_r(B_1)$.
We claim that $T \leq H$.
First suppose that $\rr_0(B_2) = r'$.
Then $Y \cdot T$ is nilpotent by Fact \ref{nilpotencepre2}.
By Fact \ref{nildecomp}, $Y$ centralizes $T$.
So $T \leq N^\o_G(Y) \leq B_2$ by Lemma \ref{Aj_Y},
 and $T \leq H$.
Next, suppose that $\rr_0(B_2) > r'$.
Then $U_0(B_2) \leq Z(H)$ by Lemma \ref{Aj_UrFB2}.
By Fact \ref{nilpotencepre2},
 $U := U_0(B_2) \cdot T$ is nilpotent.
If $r \neq \rr_0(B_2)$, then
 $T \leq C^\o_G(U_0(B_2)) \leq B_2$ by Fact \ref{nildecomp},
 and $T \leq H$.
So we may assume that $r = \rr_0(B_2)$.
If $T \not\leq B_2$, then 
 $U_{0,r}(N^\o_U(U_0(B_2))) > U_0(B_2)$ by Fact \ref{Unormalizer},
but $N^\o_G(U_0(B_2)) = B_2$, a contradiction.
Thus $T \leq H$.

Since $T \leq H$, and $U_{0,r}(H') = 1$ by Lemma \ref{Aj_UrX},
 $T$ is contained in a Carter subgroup of $H$,
 by Fact \ref{USylow_decomp}.
Now $T \leq Q$ because $T \normal H$.
Clearly $T \leq F^\o(H)$ too.
By Theorem \ref{Aj_FH} and Fact \ref{Carter_cover}-4,
 $H = Q F^\o(H) \leq C^\o_G(T)$, and hence $T \leq Z(H)$.

Now consider the case where $r > r'$.
Then $F_{r'}(B_2) \cdot T$ is nilpotent
 by Fact \ref{nilpotencepre2}.
So $F_{r'}(B_2) \leq C^\o_G(T)$ by Fact \ref{nildecomp},
 a contradiction to $N^\o_G(T) = B_2$ if $T \neq 1$.

Finally consider the case where $r < r'$.
Since $T \leq Z(H)$,
 $T Y$ is abelian by Theorem \ref{Aj_FH}.
Recall that $Y \leq F(B_2)$ by Lemma \ref{Aj_Yident2}.
Let $P := U_{0,r'}(N_{F(B_2)}(Y))$.
Then $[x,h] \in Y$ for any $x\in X$ and any $h\in P$,
 and hence $[x,h] = [x,h]^t = [x,h^t]$ for any $t\in T$.
So $[h,t] = h^{-1} h^t \in C_G(X)$.
Now $[P,T] \leq Y$ by Lemma \ref{Aj_NX} and Fact \ref{Unilcommutator}. 
Since $P$ is nilpotent, and $T$ commutes with $Y$,
 the product $T P$ is nilpotent.
By Fact \ref{nildecomp}, $P \leq N_G^\o(T) = B_1$ if $T \neq 1$,
 in contradiction with Lemma \ref{Aj_Y}. 
\end{proof}

As a result, $r'$ is also uniquely determined by $B_1$.

\begin{corollary}\label{Aj_B1rprime}
$r'$ is the minimal reduced rank in $F(B_1)$.
\end{corollary}


\begin{corollary}\label{Aj_SylowUir}
For $r \leq \rr_0(B_2)$, a Sylow $U_{0,r}$-subgroup of $H$
 is a Sylow $U_{0,r}$-subgroup of $B_1$.
\end{corollary}

\begin{proof}
By Lemma \ref{Aj_FrB1} and \ref{Aj_BNX},
 $F_r(B_1) \leq H$.
Since $Q$ is a Carter subgroup of $B_1$ by Theorem \ref{Aj_Carter},
$U_{0,r}(Q) F_r(B_1) \leq H$ is a Sylow $U_{0,r}$-subgroup
 of $B_1$ by Fact \ref{USylow_decomp}.
\end{proof}

%

\subsection{Nonabelian intersections}

In closing, we can build upon Corollary \ref{Aj_FB1CX}
 to produce a characterization of $B_1$.
In \S\ref{sec:Conclusions},
this fact will be used, together with Proposition \ref{Aj_B1B2unique},
to show that all reasonable notions of maximal intersection are equivalent.

\begin{lemma}\label{Aj_B1unique}
The subgroup $C^\o_G(X)$ is nonnilpotent.
If $H$ is nonabelian, then 
$B_1$ is the only Borel subgroup containing $C^\o_G(X)$.
\end{lemma}

\begin{proof}
By Lemma \ref{Aj_NX}, $C^\o_G(X) \leq B_1$.
By Lemma \ref{Aj_Yident} and Theorem \ref{Aj_FH},
 $U_{0,r'}(Q) \leq C^\o_G(X)$.
By Fact \ref{nildecomp}, $U_0(B_1) \leq C^\o_G(X)$ too.
By Theorem \ref{Aj_Carter} and Fact \ref{nildecomp},
 $U_{0,r'}(Q) \cdot U_0(B_1)$ is nonnilpotent.
So $C^\o_G(X)$ is nonnilpotent.

We now suppose that $H$ is nonabelian.
Suppose also that another Borel subgroup of $G$, distinct from $B_1$,
 contains $C^\o_G(X)$.
So there is a maximal pair $B_3,B_4$ whose intersection contains $C^\o_G(X)$.
We may assume that $\rr_0(B_3) \geq \rr_0(B_4)$.
Let $K := C^\o_G(X)'$.
By Corollary \ref{Aj_FB1CX}, $F^\o(B_1) \leq C^\o_G(X)$.
So $C^\o_G(X) \normal B_1$ by Fact \ref{fittingquotient},
Now $N^\o_G(K) = B_1$.
Since $K \leq F(B_3) \cap F(B_4)$,
 we have $F^\o(B_3) \leq C^\o_G(K)$,
 by Corollary \ref{Aj_FB1CX} (applied to the pair $B_3,B_4$).
Thus $\rr_0(B_1) \geq \rr_0(B_3) > \rr_0(B_4)$ by Lemma \ref{Aj_B1}.
But, as $F^\o(B_1) \leq B_4$,
 $\rr_0(B_4) \geq \rr_0(B_1)$ by Theorem \ref{nilpotence},
 a contradiction.
\end{proof}

\begin{corollary}\label{Aj_NXplus}
Suppose that $H$ is not abelian.  Then, for nontrivial definable
$X_1 \leq X$, $B_1$ is the only Borel subgroup containing $C^\o_G(X_1)$.
\end{corollary}

\begin{lemma}\label{Aj_B2unique}
Let $B$ be a Borel subgroup of $G$, distinct from $B_1$.
Suppose that $H_1 := (B \cap B_1)^\o$ is nonabelian,
 $B_1,B$ is a maximal pair, and $\rr_0(B_1) \geq \rr_0(B)$.
Then $B$ is $F^\o(B_1)$-conjugate to $B_2$.
\end{lemma}

\begin{proof}
We will apply the preceding results with the maximal pair $B_1,B$.
We observe that $H_1' \leq F(B_1) \cap F(B)$.
By Corollary \ref{Aj_B1rprime} and Theorem \ref{Aj_UrX},
 $r' = \rr_0(H_1')$, and $H',H_1' \leq F_{r'}(B_1)$. 
By Lemma \ref{Aj_BNX}, $F_{r'}(B_1) \leq H,H_1$.
Let $Q$ and $Q_1$ be Carter subgroups of $H$ and $H_1$, respectively.
By Lemma \ref{Aj_Carter}.
 $Q$ and $Q_1$ are Carter subgroups of $B_1$.
By Fact \ref{Carter_conj}-3,
 $Q^h = Q_1$ for some $h\in B_1$.
We may assume $h\in F^\o(B_1)$ by Fact \ref{Carter_cover}-4.
By Fact \ref{Carter_cover}-4 and Lemma \ref{Aj_BNX},
 $H^h = F_{r'}(B_1) Q^h = F_{r'}(B_1) Q_1 = H_1$.
Since $H_1$ is nonabelian, $B_2^h = B$ by Proposition \ref{Aj_B1B2unique}.
\end{proof}

\section{Conclusions}\label{sec:Conclusions}

Our main task is to understand how the results of \S\ref{sec:Aj}
translate down to nonmaximal intersections.  This translation is
immediate for results which describe the internal structure of $H$,
such as Theorems \ref{Aj_UrX} (applied to $H'$, not $X$) and \ref{Aj_FH},
but such direct translations are not possible for a number of results,
such as Lemmas \ref{Aj_B1} and \ref{Aj_UrFB2}.
The following summarizes the most useful consequences for
 arbitrary intersections for Borel subgroups.

\begin{proposition}\label{AjaligotC}
Let $G$ be a minimal connected simple group of finite Morley rank,
 let $B_1,B_2$ be two distinct Borel subgroups of $G$, and
 let $H$ be a connected subgroup of the intersection $B_1 \cap B_2$.
Then the following hold.
\begin{conclusions}
\item $H'$ is rank homogeneous for $r' := \rr_0(H')$, or trivial.
\item Every connected nilpotent subgroup of $H$ is abelian,.
\item $F_{r'}(H) = U_{0,r'}(H)$ is the unique Sylow $U_{0,r'}$-subgroup of $H$.
\item $F_r(H) \leq Z(H)$ for $r \neq r'$,
\item $\rr_0(C_G(H')) > \rr_0(H)$ if $H$ is nonabelian.
\end{conclusions}
\end{proposition}

\begin{proof}
We may assume $H$ is nonabelian because
 all five statements are trivial if $H$ is abelian.
Let $B_3,B_4$ be a maximal pair, containing $H$,
 with $\rr_0(B_3) \geq \rr_0(B_4)$.
The first two conclusions follow immediately from
 Theorems \ref{Aj_UrX} and \ref{Aj_FH}.
The third conclusion follows from Lemma \ref{Aj_Yident}.
For the fourth conclusion,
 $F_r(H)$ lies in a Carter subgroup $Q$ by Fact \ref{USylow_decomp},
 and $H \leq Q H' \leq C^\o_H(F_r(H))$ by Fact \ref{Carter_cover}-4.
By Lemma \ref{Aj_B1}, $\rr_0(B_3) > \rr_0(B_4) \geq \rr_0(H)$.
By Fact \ref{nildecomp}, $U_0(B_3) \leq C_G(H')$,
 and conclusion five follows.
\end{proof}

\begin{corollary}\label{AjaligotCcor}
Let $G$ be a minimal connected simple group of finite Morley rank.
Then a definable connected nonabelian nilpotent subgroup $H$ of $G$
is contained in exactly one Borel subgroup of $G$.
In particular, a nonabelian Carter subgroup of any Borel subgroup of $G$
is a Carter subgroup of $G$ itself. 
\end{corollary}

In \S\ref{sec:Aj}, Lemmas \ref{Aj_B1} and \ref{Aj_UrFB2}
told us much about the Borel subgroups involved, and this information
is lost in Proposition \ref{AjaligotC}-(4,5).
Instead, we prove that, in the nonabelian case, all
reasonable notions of maximal intersections are equivalent.
In practice, this allows one to make direct use of the analysis of \S\ref{sec:Aj}.

\begin{theorem}\label{maximalities}
Let $G$ be a minimal connected simple group of finite Morley rank,
 and let $B_1,B_2$ be two distinct Borel subgroups of $G$.
Suppose that $H := (B_1 \cap B_2)^\o$ is nonabelian.
Then the following are equivalent.
\begin{conclusions}
\item $B_1$ and $B_2$ are the only Borel subgroups of $G$ containing $H$.
\item If $B_3$ and $B_4$ are distinct Borel subgroups of $G$ containing $H$,
  then\\ $(B_3 \cap B_4)^\o = H$.
\item If $B_3 \neq B_1$ is a Borel subgroup containing $H$,
  then $(B_3 \cap B_1)^\o = H$.
\item $C^\o_G(H')$ is contained in $B_1$ or $B_2$.
\item $B_1$ and $B_2$ are not conjugate under $C^\o_G(H')$.
\item $\rr_0(B_1) \neq \rr_0(B_2)$.
\end{conclusions}
\end{theorem}

The second and third clauses express the maximality of $H$ in two different
senses.  The second clause corresponds to the sense of maximality used
in \S\ref{sec:Aj}, while the third clause has an a priori weaker sense.
The first clause goes far beyond maximality, to assure uniqueness of the
Borel subgroups, and the fourth and fifth clauses provide means to recognize
maximal intersections ``in the wild.''  We observe that the first and second
clauses are equivalent, thanks to Proposition \ref{Aj_B1B2unique}.
To prove this theorem, we first treat the most subtle implication ($4 \implies 1$)
with a lemma.

\begin{lemma}\label{maximalities_L}
Let $G$ be a minimal connected simple group of finite Morley rank,
 and let $B_1,B_2$ be two distinct Borel subgroups of $G$.
Suppose that the intersection $H := (B_1 \cap B_2)^\o$ is nonabelian,
 and that $C^\o_G(H')$ is contained in $B_1$.
Then $B_1$ and $B_2$ are the only Borel subgroups containing $H$.
\end{lemma}

\begin{proof}
Suppose toward a contradiction that $B$ is a Borel subgroup of $G$
 containing $H$ which is distinct from $B_1$ and $B_2$.
We may choose $B$ such that $H_2 := (B_2 \cap B)^\o$ is maximal,
 subject to $B \neq B_1,B_2$ and $B \geq H$.
Consider a maximal pair $B_3,B_4$ containing $H_2$,
 with $\rr_0(B_3) \geq \rr_0(B_4)$.
By Corollary \ref{Aj_NXplus},
 $B_3$ is the only Borel subgroup containing $C^\o_G(H')$,
so $B_1 = B_3$.  Thus $H = H_2$.
By Proposition \ref{Aj_B1B2unique},
 $B_1 = B_3$ and $B_4$ are the only Borel subgroups containing their intersection.
So we may assume that $B_4 \neq B_2$, as otherwise we are done.
Therefore we may also assume that $B = B_4$.
So $H_1 := (B_1 \cap B)^\o = (B_3 \cap B_4)^\o$
 is a maximal intersection containing $H_2$,
and we are free to apply \S\ref{sec:Aj} here.
We observe that $r' := \rr_0(H') = \rr_0(F(B_1) \cap F(B))$,
 by Theorem \ref{Aj_UrX}.

We first consider the case where $F_{r'}(B_2) \leq B_1$.
Since $H'$ is rank homogeneous, Fact \ref{nildecomp} says
 $F_r(B_2) \leq C^\o_G(H') \leq B_1$ for $r \neq r'$, and
 $\Ftor^\o(B_2) \leq C^\o_G(H') \leq B_1$.
So $F^\o(B_2) \leq B_1$.
Hence $F^\o(B_2) \leq H$, and $H \normal B_2$.
By Corollary \ref{Aj_FB1CX} (for $B_1,B$),
 $F^\o(B_1) \leq N^\o_G(H') = B_2$.
So $U_0(B_1) \leq H$ by Fact \ref{nilpotence},
 and $\rr_0(B_1) \leq \rr_0(H_1)$.
But this contradicts Lemma \ref{Aj_B1}.

We next consider the case where $F_{r'}(B_2) \not\leq B_1$.
Let $P := F_{r'}(H)$, $Y := F_{r'}(H_1)$, and $M := N^\o_G(P)$.
By Proposition \ref{AjaligotC}-3,
 $P = U_{0,r'}(H)$ is a Sylow $U_{0,r'}$-subgroup of $H$.
Since $P$ normalizes $F_{r'}(B_2)$,
 $P \cdot F_{r'}(B_2)$ is nilpotent by Fact \ref{nilpotencepre2}.
By Fact \ref{Unormalizer},
 $U_{0,r'}(N^\o_{P \cdot F_{r'}(B_2)}(P)) > P$. 
Since $P = U_{0,r'}(H)$, $(M \cap B_2)^\o > H$.
Since $H = (B_1 \cap B_2)^\o$, $(M \cap B_2)^\o \not\leq B_1$.
So $M$ is contained in $B_2$, by maximality of $H_2$ ($= H$).

By Lemma \ref{Aj_UFB2}, $F_{r'}(B) \not\leq H_1$.
Since $P$ normalizes $F_{r'}(B)$,
 $P \cdot F_{r'}(B)$ is nilpotent by Fact \ref{nilpotencepre2}.
By Fact \ref{Unormalizer},
 $U_{0,r'}(N^\o_{P \cdot F_{r'}(B)}(P)) > P$ too.
As $M \leq B_2$, we have
 $U_{0,r'}(N^\o_{P \cdot F_{r'}(B)}(P)) \leq H_2 = H$,
 contradicting the fact that $P$ is a Sylow $U_{0,r'}$-subgroup of $H$.
\end{proof}

\smallskip

\begin{proof}[Proof of Theorem \ref{maximalities}]
The first three clauses are successively weaker,
 and clause 4 implies clause 5.
By Lemmas \ref{Aj_B1} and \ref{Aj_B2},
 clause 2 implies clause 6.
Clearly, clause 6 implies clause 5.
So we shall examine the implications 3 to 4, and 5 to 1.
Let $B_c$ be a Borel subgroup containing $N^\o_G(H')$.

We first assume clause 3, and show clause 4.
Let $B_x$ denote $B_1$, unless $B_c = B_1$,
 in which case we let $B_x$ denote $B_2$.
By clause 3, $H := (B_c \cap B_x)^\o$.
By Lemma \ref{maximalities_L} (for $B_c,B_x$),
 $B_c \geq C^\o_G(H')$ must be one of $B_1$ or $B_2$,
 so clause 4 holds.

We now assume clause 1 fails, and show that clause 5 fails.
Then, for $i=1,2$, $C^\o_G(H') \not\leq B_i$,
 by Lemma \ref{maximalities_L}.
But $B_c,B_1$ and $B_c,B_2$ are maximal pairs,
 by Lemma \ref{maximalities_L} again.
So $\rr_0(B_c) \geq \rr_0(B_1), \rr_0(B_2)$,
 by Lemma \ref{Aj_B1unique}.
By Lemma \ref{Aj_B2unique},
 $B_1$ is $F^\o(B_c)$-conjugate to $B_2$.
By Corollary  \ref{Aj_FB1CX},
 $F^\o(B_c) \leq C^\o_G(H')$, as desired.
\end{proof}

We can summarize \S\ref{sec:Aj}, in the nonabelian case, as follows.

\begin{theorem}\label{summary}
Let $G$, $B_1$, $B_2$, and $H$ satisfy the hypotheses and conditions
of Theorem \ref{maximalities}, and let $r' = \rr_0(H')$.
Suppose that $\rr_0(B_1) \geq \rr_0(B_2)$.
Then the following hold.
\begin{conclusions}
\item $\rr_0(B_1) > \rr_0(H) = \rr_0(B_2)$, and $N^\o_G(H) = H$.
\item Every connected nilpotent subgroup of $H$ is abelian.
\item $F_{r'}(H) = U_{0,r'}(H)$ is the unique Sylow $U_{0,r'}$-subgroup of $H$.
  It is contained in $F(B_2)$, and its normalizer\o is contained in $B_2$.
\item $F_r(B_2) \leq Z(H)$ for $r \neq r'$,
  and $F_{r'}(B_2)$ is nonabelian.
  So $F_r(B_2)$ is nonabelian iff $r = r'$.
\item Carter subgroups of $H$ are Carter subgroups of $B_1$,
  and their normalizers\o are contained in $B_2$.
\item $F_{r'}(B_1) = F(B_1) \cap F(B_2)$ is rank homogeneous for $r'$,
  and $B_1$ is the only Borel subgroup containing $C^\o_G(F_{r'}(B_1))$.
\item $F_r(B_1) = 1$ for $r \leq \rr_0(B_2)$ with $r \neq r'$.
  So $r'$ is the minimal reduced rank in $F(B_1)$.
\item $F^\o(B_1)$ and $F^\o(B_2)$ are divisible.
\end{conclusions}
\end{theorem}

\begin{proof}
Part 1 follows from Lemmas \ref{Aj_B1}, \ref{Aj_B2}, and \ref{Aj_NH}.
Part 2 is Theorem \ref{Aj_FH}.
Part 3 consists of Lemmas \ref{Aj_Yident} and \ref{Aj_Y}.
Part 4 includes Lemmas \ref{Aj_UrFB2} and \ref{Aj_UFB2},
 along with Corollary \ref{Aj_B2rprime}.
Part 5 restates Theorem \ref{Aj_Carter}.
Part 6 summarizes Lemma \ref{Aj_BNX}, Theorem \ref{Aj_UrX},
 and Lemma \ref{Aj_B1unique}.
Part 7 is Lemma \ref{Aj_FrB1}, along with Corollary \ref{Aj_B1rprime}.
Part 8 is Lemmas \ref{Aj_B2div} and \ref{Aj_B1div}.
\end{proof}

\section{Genericity of Carter subgroups}\label{sec:Future}

We conclude this article by discussing an important open question:
whether a group $G$ of finite Morley rank possesses a Carter subgroup
whose conjugates are generic in $G$.  It is natural to ask whether the
structure imposed by a nonabelian intersection can be used to prove
the genericity of any related Carter subgroup.  We show that genericity
of conjugates holds for a Carter subgroups of a nonabelian intersection
iff it is a Carter subgroup of both sides.

\begin{theorem}\label{Aj_genericity}
Let $G$ be a minimal simple group of finite Morley rank,
and let $B_1,B_2$ be a maximal pair of Borel subgroups of $G$.
Suppose that $H := (B_1 \cap B_2)^\o$ is nonabelian,
 and that $\rr_0(B_1) > \rr_0(B_2)$.
If $Q$ is a Carter subgroup of $H$, then the following are equivalent.
\begin{conclusions}
\item $Q$ is a Carter subgroup of $B_2$.
\item $Q$ is a Carter subgroup of $G$.
\item $\bigcup Q^G$ is generic in $G$.
\item $\bigcup H^G$ is generic in $G$.
\end{conclusions}
\end{theorem}

To prove this, we recall the following fact.

\begin{fact}[{cf. \cite[\qLemma 3.3]{CJ01}}]\label{f:BGgen_L}
Let $G$ be a connected group of finite Morley rank, and
 let $C$ be a definable almost self-normalizing subgroup.
 Suppose there is a definable subset $J$ of $C$, not generic in $C$,
 such that $C \cap C^g \subseteq J$ whenever $g\notin N_G(C)$.
Then $\cup (C \setminus J)^G$ is generic in $G$.
\end{fact}

We also need a lemma due to \Frecon.

\begin{fact}[{\Frecon, cf.\ \cite[\qLemma 3.5]{CJ01}}]\label{f:QgenB_L}
Let $B$ be a connected solvable group of finite Morley rank,
 and let $C$ be a Carter subgroup of $B$.
Then there is a definable subset $J$ of $C$, not generic in $C$,
 such that $C \cap C^b \subseteq J$ whenever $b\notin N_B(C)$.
\end{fact}

These two facts can be used to prove the following.

\begin{lemma}\label{Aj_genericity_L}
Let $G$ be a minimal connected simple group of finite Morley rank.
Let $H$ be a connected solvable subgroup of $H$, and
 let $Q$ be a Carter subgroup of $H$.
Then $\bigcup H^G$ is generic in $G$ iff $\bigcup Q^G$ is generic in $G$.
\end{lemma}

\begin{proof}
The if direction is immediate because $\rk(\bigcup H^G) \geq \rk(\bigcup Q^G)$.
So we assume $\rk(\bigcup Q^G) < \rk(G)$,
 and show that $\rk(\bigcup H^G) < \rk(G)$.
Let $H_* := H \setminus \bigcup Q^H$.
By Fact \ref{f:QgenB_L}
 there is a set $J \subset Q$, not generic in $Q$,
 such that $Q \cap Q^h \subset J$ for $h\in H$.
By Fact \ref{f:BGgen_L},
 $\bigcup Q^H$ is generic in $H$, and $\rk(H_*) < \rk(H)$.
So $\rk(\bigcup H_*^G) \leq \rk(G/N_G(H) \times H_*)
 = \rk(G) - \rk(N_G(H) + \rk(H_*) < \rk(G)$.
Since $H = H_* \sqcup \bigcup Q^G$, we have $\rk(\bigcup H^G)
 \leq \max(\rk(\bigcup H_*^G),\rk(\bigcup Q^G)) < \rk(G)$, as desired.
\end{proof}

We will also need the following two tools.

\begin{fact}\label{decent_tori}
Let $G$ be a minimal simple group of finite Morley rank,
 and let $Q$ be a Carter subgroup of $G$ with divisible $p$-torsion.
Then $\bigcup Q^G$ is generic in $G$.
\end{fact}

\begin{proof}
Let $P$ be a divisible abelian $p$-subgroup of $Q$.
Then $P$ is contained in a maximal decent torus $T$ (see \cite{Ch05}).
By Fact \ref{Carter_Sylow}-5,
 $N^\o_G(T)$ contains a Carter subgroup $Q_1$ of $C^\o_G(P)$.
Since $Q$ is a Carter subgroup of $C^\o_G(P)$ too,
 we may assume $Q = Q_1$ by Fact \ref{Carter_conj}-3.
So $T \leq Q$ by Fact \ref{cent_divtorsion},
 and $Q$ is a Carter subgroup of $C^\o_G(T)$.
By \cite{Ch05}, $\bigcup C^\o_G(T)^G$ is generic in $G$,
So $\bigcup Q^G$ is generic in $G$, by Lemma \ref{Aj_genericity_L}.
\end{proof}

\begin{fact}\label{nildecomp_derived}
Let $Q$ be a nilpotent group of finite Morley rank.  Then
$$ Q' = U_{0,1}(Q)' \cdot U_{0,2}(Q)' \cdots U_{0,\rk(Q)}(Q)' \cdot U_2(Q)' \cdot U_3(Q)' \cdots $$
\end{fact}

\begin{proof}
Let $P$ denote the product on the right hand side.
Clearly $P \leq Q'$.  By Fact \ref{nildecomp},
 $Q/P$ is abelian, so $Q' \leq P$.
\end{proof}

\begin{lemma}\label{find_Uir_elements}
Let $Q$ be a nilpotent group of finite Morley rank without divisible torsion,
and suppose that $\Uir(Q) \neq 1$.  Then there is a generic subset $Q^*_r$
of $Q$ such that $\Uir(d(k)) \neq 1$ for all $k\in Q^*_r$. 
\end{lemma}

\begin{proof}
There is a (not necessarily unique) maximal connected normal subgroup $P$
 of $Q$ such that $\Uir(Q/P) \neq 1$.
Then $P$ contains $U_{0,r}(Q)'$ for any $r$, and $U_p(Q)'$ for any prime $p$.
So $P$ contains $Q'$ by Fact \ref{nildecomp_derived}.
Hence the normality assumption on $P$ is superfluous.
By Fact \ref{Upnilpotence}, $U_p(Q) \leq P$ for any prime $p$.
By Fact \ref{Sylow_conplus}, $Q/P$ is torsion-free.
If $\Uir(d(x P/P)) = 1$ for some $x\in Q^*_r$,
 then $\Uir((Q/P) / d(x P/P)) \neq 1$ by Lemma \ref{Uhom}.
So $\Uir(d(x P/P)) \neq 1$ by the maximality of $P$.
By Fact \ref{Uhom}, $\Uir(d(x)) \neq 1$, as desired.
\end{proof}

\smallskip

\begin{proof}[Proof of Theorem \ref{Aj_genericity}]
We first observe that $B_1$ and $B_2$ satisfy Hypothesis \ref{Aj_hypothesis}.
So clause 1 and clause 2 are equivalent, by Lemma \ref{Aj_Y}.
Lemma \ref{Aj_genericity_L} proves that clauses 3 and 4 are equivalent.

We next assume clause 2 and prove clause 3.
Suppose toward a contradiction that $\bigcup Q^G$ is not generic in $G$.
Then $Q$ contains no divisible torsion, by Fact \ref{decent_tori}.
By Lemma \ref{Aj_Carter}, $U_{0,r'}(Q) \neq 1$.
So there is a generic subset $Q^*_{r'}$ of $Q$ such that
 $U_{0,r'}(d(h)) \neq 1$ for any $h\in Q^*_{r'}$,
 by Lemma \ref{find_Uir_elements}. 
For $i=1,2$, there is a set $J_i \subset Q$, not generic in $Q$,
 such that $Q \cap Q^b \subset J_i$ for $b\in B_i$,
 by Fact \ref{f:QgenB_L}.
Then $Q_* := Q^*_{r'} \setminus (J_1 \cup J_2)$ is a generic subset of $Q$.
By Fact \ref{f:BGgen_L},
 there is an $h \in Q_* \cap Q^g$, for some $g\notin N_G(Q)$.
Since $h \in Q^*_{r'}$, $K := U_{0,r'}(d(h)) \neq 1$.
By Theorem \ref{Aj_FH} and Lemma \ref{Aj_BNX},
 $H',(H^g)' \leq C^\o_G(K)$.
Since $h\in Q \cap Q^g$, and $Q$ is abelian,
 we have $Q,Q^g \leq C^\o_G(K)$ too,
 so $H,H^g \leq C^\o_G(K)$.
So $Q^g \leq H^g$ lies in either $B_1$ or $B_2$
 by Proposition \ref{Aj_B1B2unique}.
Since $Q^g$ is a Carter subgroup of either $B_1$ or $B_2$,
 $h \in Q \cap Q^g$ lies in either $J_1$ or $J_2$ by Fact \ref{Carter_conj}-3,
 a contradiction.

We next assume clause 2 fails, and show that clause 3 fails.
For any group $Q$,  $\rk(\bigcup_{g\in G} Q) \leq
 \rk(G/N_G(Q) \times Q) = \rk(G) - \rk(N_G(Q)) + \rk(Q)$.
So $\rk(G) - \rk(\bigcup_{g\in G} Q) \geq \rk(N_G(Q)/Q) > 0$,
 as desired.
\end{proof}

\begin{corollary}\label{Aj_genericity_cor}
Let $G$ be a minimal connected simple group of finite Morley rank,
 and let $B_1,B_2$ be two distinct Borel subgroups of $G$.
Suppose that $H := (B_1 \cap B_2)^\o$ is nonabelian, and
 contains a Carter subgroup $Q$ of both $B_1$ and $B_2$.
Then $\rr_0(B_1) \neq \rr_0(B_2)$,
 and $\bigcup Q^G$ is generic in $G$.
\end{corollary}

\begin{proof}
We may assume that $\rr_0(B_1) \geq \rr_0(B_2)$.
If $\rr_0(B_1) \neq \rr_0(B_2)$,  then
 $B_1$ and $B_2$ are the only Borel subgroups containing $H$,
 by Theorem \ref{maximalities},
and hence Theorem \ref{Aj_genericity} applies.

Suppose toward a contradiction that $\rr_0(B_1) = \rr_0(B_2)$.
Let $B_3$ be a Borel subgroup containing $N^\o_G(H')$.
For $i=1,2$,
 $B_i$ and $B_3$ are the only Borel subgroups containing
 their intersection $H_i := (B_i \cap B_3)^\o$, by Theorem \ref{maximalities}.
Since $Q$ is a Carter subgroup of $B_i$, and $Q \leq H \leq B_3$,
 $Q$ is a Carter subgroup of $H_i$ too.
By Theorem \ref{Aj_Carter}, $Q$ is a Carter subgroup of $B_3$.
By Lemma \ref{Aj_BNX}, $H_i' \leq F_{r'}(B_3) \leq H_i$.
By Fact \ref{Carter_cover}-4,
 $H_1 = F_{r'}(B_3) Q = H_2$, contradicting $B_1 \neq B_2$.
\end{proof}

Even if $Q$ fails to be a Carter subgroup of $B_2$, one might hope to
show that the conjugates of a Carter subgroup $C$ of $B_2$ are generic.
However, a nonabelian intersection appears to place few
constraints on $C$.

On the positive side,
Corollary \ref{AjaligotCcor} and Fact \ref{Carter_cleanup}-2
show that a nonabelian Carter subgroup of any Borel subgroup of
a minimal connected simple group $G$ of finite Morley rank is
actually a Carter subgroup of $G$ itself.  So it is natural to ask
when these Carter subgroups are generic in $G$.

\begin{proposition}\label{Carter_noneverything}
Let $G$ be a minimal simple group of finite Morley rank,
 and let $Q$ be a nonabelian Carter subgroup of $G$.
If $Q$ is not rank homogeneous, then $\bigcup Q^G$ is generic in $G$.
\end{proposition}

We observe that nilpotent Borel subgroups with unipotent torsion are generic.

\begin{fact}\label{generic_nilpotent_Borel}
Let $G$ be a minimal simple group of finite Morley rank
 with a nilpotent Borel subgroup $B$ which is not divisible.
Then $\bigcup B^G$ is generic in $G$.
\end{fact}

\begin{proof}
By Lemma \ref{Upjaligot}, $B \cap B^g = 1$ for $g\notin N_G(B)$.
So genericity follows by Fact \ref{f:BGgen_L}.
\end{proof}

\smallskip

\begin{proof}[Proof of Proposition \ref{Carter_noneverything}]
Suppose toward a contradiction that $\bigcup Q^G$ is not generic in $G$,
 but that $Q$ is not rank homogeneous.
Since $Q$ is nonabelian, there is a unique Borel subgroup $B$
 containing $Q$, by Corollary \ref{AjaligotCcor}.
By Fact \ref{f:QgenB_L},
 there is a set $J_0 \subset Q$, not generic in $Q$, such that
 $Q \cap Q^b \subset J_0$ for any $b\in B \setminus N_G(Q)$.

By Fact \ref{decent_tori}, $Q$ contains no divisible torsion.
If $Q$ had bounded exponent, $Q$ would be a Borel subgroup,
 by Fact \ref{Upnilpotence},
in contradiction with Fact \ref{generic_nilpotent_Borel}.
So $Q$ is not of bounded exponent.

Since $Q$ is nonabelian, either $U_{0,r}(Q)' \neq 1$ for some $r$,
 or $U_p(Q)' \neq 1$ for some $p$ prime, by Fact \ref{nildecomp}.
We take $p=0$ in the former case and $r = 0$ in the latter.
First consider the case where $U_{0,s}(Q) \neq 1$ for some $s \neq r$.
Then there is a generic subset $Q^*_s$ of $Q$ such that
 $U_{0,s}(d(h)) \neq 1$ for all $h\in Q^*_s$,
 by Lemma \ref{find_Uir_elements}. 
Next consider the case where $U_{0,s}(Q) = 1$ for all $s \neq r$.
We may assume that $p=0$,
 since $Q$ is not of bounded exponent.
Since $Q$ is not rank homogeneous,
  there is a prime $q$ such that $U_q(Q) \neq 1$.
So there is a maximal connected subgroup $P$ of $Q$
 such that $U_q(Q/P) \neq 1$.
Clearly, $U_q(Q/P) = Q/P$.
So, for any $h \notin P$, $d(k)$ contains a $q$-element.
In either case, there is a generic subset $Q_*$ of $Q$ given by
 either $Q_* := Q^*_s \setminus J_0$ or $Q_* := Q \setminus (P \cup J_0)$.

By Fact \ref{f:BGgen_L},
 there is an $h \in Q_* \cap Q^g$, for some $g\notin N_G(Q)$.
Since $h\notin P$, either $K := U_{0,s}(d(h)) \neq 1$,
 or $d(h)$ contains a $q$-torsion subgroup $K$.
By Fact \ref{nildecomp},  $C^\o_B(K)$ contains
 the nonabelian nilpotent subgroup $U_{0,r}(Q)$, or $U_p(Q)$.
So $B$ is the only Borel subgroups containing $C^\o_G(K)$,
 by Proposition \ref{AjaligotC}.
But $Q,Q^g \leq C^\o_G(K)$.
So they are $B$-conjugate by Fact \ref{Carter_conj}-3,
 a contradiction to $h\notin J_0$.
\end{proof}

\section*{Acknowledgments}

I thank my advisor Gregory Cherlin for direction during my thesis work,
of which this article is an outgrowth, and for guidance during its later
developments.  I also thank Tuna \Altinel, Alexandre Borovik, Eric Jaligot,
and Olivier \Frecon, for corrections and suggestions.
The applications of this paper in \cite{BCJ} are critical to our
understanding of its material, and this would also be impossible without
the collaboration of Gregory Cherlin and Eric Jaligot on \cite{BCJ}.

Financial support for this work comes from an NSF Graduate Research
Fellowship, NSF grant DMS-0100794, DFG grant Te 242/3-1, and
the Isaac Newton Institute, Cambridge.
I also thank the following institutions for their hospitality:
University of Birmingham, University of Manchester,
IGD at Universit\'e Lyon I, Isaac Newton Institute, Cambridge,
and CIRM at Luminy.

\small
\bibliographystyle{elsart-num}
\bibliography{burdges,fMr}

\end{document}